\documentclass{siamltex}
\usepackage{amsmath,amssymb,url}
\author{Erik G. Boman\thanks{Sandia National Labs, 
Albuquerque, NM 87185-1318, USA, egboman@sandia.gov.} 
  \and Bruce Hendrickson\thanks{Sandia National Labs, 
Albuquerque, NM 87185-1318, USA, bah@cs.sandia.gov.} \and 
Stephen Vavasis\thanks{Department of Combinatorics and Optimization,
University of Waterloo, Waterloo, ON N2L 2W2, Canada, 
vavasis@math.uwaterloo.ca.  Corresponding author.
Part of this work was done while the
author was in the Department of Computer Science of Cornell University.}}
\title{Solving Elliptic Finite Element Systems in Near-Linear Time
with Support Preconditioners\thanks{The work of the first two authors 
was funded by
the Applied Mathematical Sciences program,
U.S.\ Department of Energy, Office of Energy Research and performed
at Sandia National Labs, a multiprogram laboratory
operated by Sandia Corporation, a Lockheed Martin Company,
for the U.S.\ Dept.~of Energy's National Nuclear Security Administration
under contract number DE-AC-94AL85000.
The third author was supported in part
by NSF Grant CCF-0085969 and by a grant from NSERC, Canada.}
}
\date{January 2007}

\def\norm2#1{\|#1\|_2}
\newcommand\bxi{\mbox{\boldmath{$\xi$}}}
\newcommand\bz{{\bf 0}}
\newcommand\eref[1]{$(\ref{#1})$}
\newcommand\f{{\mathbf{f}}}
\newcommand\I{{\mathbf{I}\kern-.6ex\mathbf{I}}}
\newcommand\Q{{\mathcal{Q}}}
\newcommand\R{{\mathbb{R}}}
\newcommand\T{{\mathcal T}}
\newcommand\TT{{\bf T}}
\renewcommand\v{{\bf v}}
\newcommand\x{{\bf x}}
\newcommand\y{{\bf y}}
\newtheorem{assumption}{Assumption}
\newcommand\commentout[1]{\relax}
\newcommand{\nnu}{k}

\begin{document}
\maketitle
\begin{abstract}
We consider linear systems arising
from the
use of the finite element method for solving scalar
linear elliptic problems.  Our main result is that
these linear systems, which are
symmetric and positive semidefinite,
are well approximated by symmetric diagonally dominant matrices.
Our framework for defining matrix approximation is support
theory.  Significant graph theoretic work has already
been developed in the support framework
for preconditioners in the diagonally dominant
case, and in particular it is known that
such systems can be solved with iterative methods in nearly linear time.
Thus, our approximation result implies that these graph theoretic
techniques can also solve a class of finite element problems
in nearly linear time.
We show that the support number bounds, which control
the number of iterations in the preconditioned iterative
solver, 
depend on mesh quality measures but not on
the problem size or shape of the domain.
\end{abstract}

\section{Introduction}
\label{sec:intro}
Finite element discretizations of elliptic partial differential
equations (PDEs) give rise
to large sparse linear systems of equations. A topic of
great interest is preconditioners for iterative solution 
of such systems.  We focus on scalar boundary value problems
of the form
$\nabla  \cdot(\theta \nabla u)=-f$, in which
$\theta$ is a scalar conductivity field.  See \eref{bvp}
below for a more detailed statement of the PDE under consideration.
Such PDEs
arise in a variety of physical applications listed in
Section~\ref{sec:fem}.   We prove that the stiffness matrix
$K$ of this PDE, which is defined precisely by \eref{Kijq} below,
can be well approximated by a diagonally dominant matrix $\bar K$.
``Well approximated'' means that $\kappa(\bar K,K)$ has an upper bound in terms
of mesh quality measures only; the notation $\kappa(\cdot,\cdot)$
is defined below.
Since significant theory has been developed for diagonally
dominant matrices, our result shows that the same theory extends
to this class of stiffness matrices.
In particular, our approximation result means that the
system $K\x=\f$ arising from the finite element method (FEM) can
be solved in nearly linear time by preconditioned iterative methods.

Our analysis uses the \emph{support theory} framework
described in \cite{BH04} 
for analyzing condition numbers and (generalized) eigenvalues for
preconditioned systems.  
Our analysis is as follows.  In
Sections~\ref{sec:matrixapprox}--\ref{sec:secondfac}, we state and prove
our theorem that $K$ may be factored
as
$K=A^T\bar{D}^{1/2}H\bar{D}^{1/2}A$.  A
preliminary  general-purpose result in
Section~\ref{sec:support} shows that a matrix of this form
can be approximated by $\bar K=A^T\bar{D}A$ (a diagonally dominant
matrix), with the quality of approximation depending on $\kappa(H)$.
The analysis of $\kappa(H)$ in Section~\ref{sec:secondfac} establishes
that this quantity depends on the space dimension $d$ and
degree $p$ of the FEM, the quadrature rule, and various mesh quality measures
but does not depend on the number of nodes or elements or on the
shape of the domain or size of elements.  It also does not depend on
the conductivity field $\theta$ under certain assumptions to be 
made below.

The idea of approximating FEM systems by diagonally dominant matrices
is not new; see for example Gustafsson \cite{Gus94}. In fact, our approach
is similar to Gustafsson's in that we approximate each element
matrix by a diagonally dominant matrix. In contrast to \cite{Gus94},
we are able to rigorously prove bounds on the spectral
properties of our approximation (and thus also the preconditioner). 

An approach in \cite[\S7.3]{Axelsson} proposes to precondition general
non-diagonally dominant FEM problems with a diagonally dominant preconditioner
obtained by using a lower order method.  Again, our result is a departure
from theirs
because we get precise and rigorous 
bounds on the quality of the
approximation.

In Section~\ref{sec:comptest},
we specialize the theory to four common cases of the finite element
method and report on computational tests
for these cases.  The computational tests are intended to illustrate
the size of the constants in our main theorem
and do not measure the practical performance of our solution strategy
versus competitors.  See further
remarks on the practical efficiency of our preconditioning approach in 
Section~\ref{sec:notesaddedinrevision}.

\section{Support Theory}
\label{sec:support}
Let $A,B$ be symmetric positive semidefinite (SPSD) matrices of the same 
size, and let $N(A)$ and $N(B)$ denote their null spaces.
We
define the {\em support number} of $A$ with respect to $B$ to be
\begin{equation}
\sigma(A,B)=\sup_{\x\in\R^n-N(B)}\frac{\x^TA\x}{\x^TB\x}.
\label{supnum}
\end{equation}
Note that this quantity is finite only if $N(B)\subset N(A)$.
We follow the convention established in the previous
literature of using $\sigma(\cdot,\cdot)$ to denote support numbers.
Unfortunately, $\sigma$ is also commonly used to denote singular values. 
In this paper, $\sigma$ with one argument is a singular value,
and $\sigma$  with two arguments is a support number.

The results in this section can be partially generalized to indefinite
symmetric matrices.  This generalization requires a more general
definition of support number than \eref{supnum} and requires
the use of the Symmetric Product Support Theorem from 
\cite{BH04}.   Since this paper focuses on the positive semidefinite case,
we omit this generalization.

If $A$ and $B$ are symmetric positive definite (SPD),
then $\sigma(A,B) = \lambda_{\max}(A,B)$, the 
largest generalized eigenvalue. When $B$ is a preconditioner for $A$
in the preconditioned conjugate gradient iterative method \cite{GVL},
the condition number of the preconditioned system is given
by $\sigma(A,B) \sigma(B,A)$, which we denote as $\kappa(A,B)$.
When $A, B$ are SPD, then $\kappa(A,B) = \kappa(B^{-1}A)$ where
$\kappa$ denotes the standard spectral condition number.

A principal goal of this paper is to propose the construction of 
a symmetric diagonally dominant matrix $\bar K$ that approximates
the finite element stiffness matrix $K$ in the sense that
$\kappa(K,\bar{K})$ is not too large.  The following two lemmas
and subsequent theorem define
our framework for defining $\bar{K}$ and bounding
$\kappa(K,\bar{K})$.

\begin{lemma}
\label{lem:UV1}
Suppose  $V\in\R^{n\times p}$, and suppose $H\in\R^{p\times p}$ is SPD.
Then
\[ \sigma(VHV^T,VV^T) \le \lambda_{\max}(H), \]
where $\lambda_{\max}$ denotes the largest eigenvalue.
\end{lemma}
\begin{proof}
\begin{eqnarray*}
\sigma(VHV^T,VV^T) & = &
\sup_{\x\in \R^{n}-N(V^T)}\frac{\x^TVHV^T\x}{\x^TVV^T\x} \\
&= &
\sup_{\y\in R(V^T)}\frac{\y^TH\y}{\y^T\y} \\
&\le & \lambda_{\max}(H).
\end{eqnarray*}
The second line follows from the first by substituting $\y=V^T\x$
and the third from the Courant-Fischer minimax theorem
\cite[\S8.1.1]{GVL}.
Here, $R(V^T)$ denotes the range-space of $V^T$.
\end{proof}

\begin{lemma}
\label{lem:UV2}
Suppose  $V\in\R^{n\times p}$, and suppose $H\in\R^{p\times p}$ is SPD.
Then
\[ \sigma(VV^T,VHV^T) \le 1/\lambda_{\min}(H), \]
where $\lambda_{\min}$ denotes the smallest eigenvalue.
\end{lemma}
\begin{proof}
Using the same idea in the previous proof leads to a supremum over the
Rayleigh quotient $\y^T\y/(\y^TH\y)$, which is bounded above by
$1/\lambda_{\min}(H)$.
\end{proof}

Combining these two lemmas yields a result for the
condition number. 
\begin{theorem}
\label{thm:condno}
Suppose  $V\in\R^{n\times p}$, and suppose $H\in\R^{p\times p}$ is SPD.
Then
\[ \kappa(VHV^T,VV^T) \le \kappa(H). \]
\end{theorem}
\begin{proof}
The result follows from Lemma~\ref{lem:UV1} and \ref{lem:UV2} 
and by using the fact $\kappa(H) = \lambda_{\max}(H)/\lambda_{\min}(H)$
when $H$ is SPD.
\end{proof}

The way we will apply this theorem is to let $V=A^T\bar{D}^{1/2}$,
where $A$ is the node-arc incidence matrix of a graph and $\bar D$ is
a diagonal weight matrix.
Then $\bar{K}=VV^T=A^T\bar{D}A$ is diagonally dominant while 
$K=VHV^T=A^T\bar{D}^{1/2}H\bar{D}^{1/2}A$ is not
(in general). Therefore, we now have a tool to approximate
non-diagonally-dominant matrices.

\section{Finite element analysis}
\label{sec:fem}
In this section we provide a brief summary of isoparametric finite element
approximation as well as an introduction to notation that is crucial
for our main theorem.  The material in this section is 
standard \cite{Johnson} in textbooks, except that our description herein
uses notation  for indexing that is more detailed than usual.  This section 
concludes with a summary of our notation.

The class of problems under consideration consists
of finite-element discretizations
of the following second-order elliptic boundary value problem.  
Find $u:\Omega\rightarrow\R$ satisfying
\begin{equation}
\begin{array}{rcll}
\nabla\cdot(\theta \nabla u) &=& -f & \mbox{on $\Omega$}, \\
u &=& u_0 & \mbox{on $\Gamma_1$}, \\
\theta \partial u/\partial n &=& g & \mbox{on $\Gamma_2$}.
\end{array}
\label{bvp}
\end{equation}
Here, $\Omega$ is a bounded open subset of $\R^d$ (typically
$d=2$ or $d=3$), $\Gamma_1$ and $\Gamma_2$ form a partition of
$\partial\Omega$, $\theta$ is a given scalar field on $\Omega$ that
is positive-valued everywhere and is sometimes called the
{\em conductivity}, $f:\Omega\rightarrow\R$ is a given
function called the {\em forcing function}, $u_0$ is a given function
called the {\em Dirichlet boundary condition}
and $g$ is another given function called the {\em Neumann boundary
condition.}

This problem has applications to many problems in mathematical physics.
For example, $u$ can represent voltage in a conducting medium, in which
case $\theta$ stands for electrical conductivity.  The two types of boundary
conditions stand for, respectively, a boundary point
held at fixed voltage or a boundary point electrically
insulated (or with prespecified nonzero current).  
Another application is thermodynamics, in which $u$ stands
for temperature of the body, $\theta$ for thermal conductivity, $\Gamma_1$
for a boundary held at fixed temperature, and $\Gamma_2$ for a boundary
insulated (or with prespecified heat flow).  Problem \eref{bvp} is also
used to model membrane deflection and gravity.  It arises as a subproblem
in fluid flow modeling.

For clarity, we keep track of our assumptions 
by explicitly numbering them.  One
assumption has already been made:
\begin{assumption}
For all $x\in\Omega$, $\theta(x)>0$.
\label{postheta}
\end{assumption}

Without this assumption, the problem may be ill posed.

The first step in the isoparametric finite element method is to produce
a mesh of this domain.  For the remainder of the paper, we assume 
isoparametric elements are defined with respect to the usual polynomial
basis on a simplicial reference element, although the results can be
generalized to other reference elements and basis families.

In more detail, 
let $T_0$ denote the standard unit $d$-simplex: for $d=2$, this
simplex is the triangle with vertices $(0,0)$, $(1,0)$, $(0,1)$,  and
for $d=3$ this simplex is the tetrahedron with vertices $(0,0,0)$,
$(1,0,0)$, $(0,1,0)$, $(0,0,1)$.
Let $p$ denote the polynomial order of the finite element method.
In the reference
element, position an evenly spaced array of $l=(p+1)(p+2)/2$ 
{\em reference nodes}
when $d=2$ or of $l=(p+1)(p+2)(p+3)/6$ nodes when $d=3$.  These nodes
have coordinates of the form $(i/p,j/p)$ (in two dimensions) or
$(i/p,j/p,k/p)$ (in three dimensions) where $i,j$ (and, in three dimensions,
$k$) are nonnegative
integers whose sum is at most $p$.  Let these nodes be enumerated
$z_1,\ldots,z_l$.

One generates a
{\em mesh} of $\Omega$ composed of $m$ elements to be defined via
mapping functions.
For each element  $t=1,\ldots,m$, there is a {\em mapping function}
$\phi_t$ that maps $T_0$ to $\R^d$.
Let $\nabla\phi_t$ denote the derivative (Jacobian)
of $\phi_t$.
The following assumption is nearly
universal in the literature.

\begin{assumption}
Mapping function $\phi_t:T_0\rightarrow\R^d$ is injective, and
$\det(\nabla\phi_t(z))>0$ 
for all $z\in T_0$, $t=1,\ldots,m$.
\label{posdetphi}
\end{assumption}

The $t$th {\em element} is defined to be $\phi_t(T_0)$ and
hence has a shape of a curved simplex.
The function $\phi_t$ carries the $l$ reference nodes $z_1,\ldots,z_l$
to $l$ {\em real-space nodes} $\phi_t(z_1),\ldots,\phi_t(z_l)$, 
which are often just called {\em nodes}.
These may be denoted $\zeta_{t,1},\ldots,\zeta_{t,l}$.
The nodes are chosen so that the nodes on the boundary of
the $t$th element coincide with the corresponding nodes on the boundaries
of its neighbors.

We restrict  $\phi_t$ to be a polynomial of degree $p$, in which case it is
uniquely determined by the positions of its real-space nodes.  In
more detail, let $N_\mu$, $\mu=1,\ldots,l$,
be a real-valued  degree-$p$ polynomial function
on $T_0$ with the property
that 
\begin{equation}
N_\mu(z_\nu)=\left\{
\begin{array}{ll}
1 & \mbox{if $\mu=\nu$,} \\
0 & \mbox{if $\mu\ne \nu$}
\end{array}
\right\}
\mbox{ for all $\nu=1,\ldots,l$.} 
\label{deltafunc}
\end{equation}

It is not hard to write down an explicit formula for $N_\mu$ and to
prove that \eref{deltafunc} uniquely determines $N_\mu$ among degree-$p$
polynomials.  These functions $N_1,\ldots,N_l$ are called {\em shape
functions}.
For each $t$, 
$\phi_t$ necessarily has the formula
$\phi_t=\zeta_{t,1}N_1+\cdots+\zeta_{t,l}N_l$.

There are many duplicate entries in the list $\zeta_{1,1},\ldots,
\zeta_{m,l}$ because of
common nodes at the boundaries of the elements.
Let $w_1,\ldots,w_{n'}$ be a listing of the real-space nodes with all
duplicates removed.  
This renumbering is specified by an 
index mapping function, the {\em local-to-global} numbering
map, denoted $LG$ and carrying an index pair $(t,\mu)$
to an index $i\in 1,\ldots,n'$ so that $w_i\equiv \zeta_{t,\mu}$. 

Finally, the {\em mesh} $\T$ is specified by listing the
nodes $w_1,\ldots,w_{n'}$ and the local-global mapping defined by
$LG$.  From this data, one can deduce all of the $\phi_t$'s.

Let $\tilde\Omega$ be the union of the elements, which
may be slightly different from $\Omega$.
This is because  at the boundary of
$\Omega$, the elements may either protrude a bit outside
$\Omega$ or may fail to cover a small
part of the boundary.  Let the boundary of $\tilde\Omega$ be partitioned
into $\tilde\Gamma_1$ and $\tilde\Gamma_2$ in correspondence with the
partition $\Gamma_1,\Gamma_2$ of the boundary of $\Omega$.
It is necessary to transfer $\theta$ (and several other
fields) from $\Omega$ to $\tilde\Omega$ and to transfer the boundary
conditions to $\tilde\Gamma_1$ and $\tilde \Gamma_2$.
We omit the details.

The next step in the isoparametric finite element method is to define
{\em basis functions} $\pi_1,\ldots,\pi_{n'}$, which are
functions from $\tilde\Omega$ to $\R$
satisfying
\begin{equation}
\pi_i(w_k)=\left\{
\begin{array}{ll}
1 & \mbox{if $i=k$}, \\
0 & \mbox{if $i\ne k$}
\end{array}
\right\}
\mbox{for $k=1,\ldots,n'$.}
\label{deltafunc2}
\end{equation}
We require that $\pi_i$, when restricted to an element $t$,
must be of the form $N_\mu\circ \phi_t^{-1}$ where $LG(t,\mu)=i$ or
else must be identically 0 
if $i$ does not occur among $LG(t,1),\ldots,LG(t,l)$.
This requirement uniquely determines
the $\pi_i$'s.  It can be shown that $\pi_i$ is continuous but
fails to be differentiable at inter-element boundaries.

Let the global numbering be chosen so that
$w_1,\ldots,w_n$, the first $n$ real-space
nodes, are the nodes with no Dirichlet boundary condition imposed,
and let the final $n'-n$ nodes have Dirichlet boundary conditions.
Now we can define the {\em exact assembled stiffness matrix}
$K^{\rm exact}$ of the finite element method to be the $n\times n$ matrix
whose $(i,j)$ entry is given by
\begin{equation}
K^{\rm exact}(i,j)=
\int_{\tilde\Omega}\nabla \pi_i(x)\cdot\theta(x)\nabla\pi_j(x)\,dx
\label{Kij}
\end{equation}
where $\nabla$ as usual denotes the gradient.

Matrix $K^{\rm exact}$ 
is sparse, symmetric and positive semidefinite.  Symmetry is obvious;
semidefiniteness follows from a fairly straightforward argument
that we omit,
and sparsity follows because $K^{\rm exact}(i,j)$ 
is nonzero only if there is an element
$t$ that contains both nodes $w_i$ and $w_j$.

Integral \eref{Kij} is difficult to compute directly because evaluating
$\pi_i$ requires evaluation of $\phi_t^{-1}$.
Fortunately, this difficulty is avoided by breaking the
integral into a sum over elements then carrying out
the integral over the reference domain following a change of variables
as follows.
\begin{equation}
K^{\rm exact}(i,j)
=
\sum_{t=1}^m\int_{T_0}\nabla_x \pi_i(\phi_t(z))\cdot\theta(\phi_t(z))
\nabla_x\pi_j(\phi_t(z))\det(\nabla \phi_t(z))\,dz, \label{Kij2}
\end{equation}
where the notation $\nabla_x$ means derivative with respect to the 
coordinates
of element
$t$ (as opposed to derivative with respect to $z$, the coordinates of $T_0$).

The integrand of \eref{Kij2} is evaluated using the chain rule for derivatives.
Assume $w_i$ is a node of element
$t$ (else the above integral is 0).  Let $\mu$
be the index such that $LG(t,\mu)=i$, so that
$\pi_i=N_\mu\circ\phi_t^{-1}$ on element $t$.
Then
\begin{eqnarray}
\nabla_x \pi_i(\phi_t(z)) &=& \nabla_x N_\mu(z) \nonumber\\
& = &\nabla\phi_t(z)^{-T}\cdot \nabla_z N_\mu(z), \label{gradpi}
\end{eqnarray}
and similarly for $\pi_j$.  Here, $\nabla\phi_t(z)^{-T}$
denotes the transposed inverse of the $d\times d$ matrix $\nabla\phi_t(z)$,
which exists by 
Assumption~\ref{posdetphi}.
The `$\cdot$' notation in the previous
formula indicates matrix-vector multiplication.  
Gradients are regarded as length-$d$ column vectors.

We assume that the entries of $K$ are not the exact value
of the integral \eref{Kij2} but are
obtained by a quadrature rule that we now discuss.  Let
$r_1,\ldots,r_q$ be points in the interior of the reference element
$T_0$ called the {\em Gauss points}.  (As is common
practice, we use this terminology even
if the quadrature rule is not derived from Gaussian quadrature.)
Let $\omega_1,\ldots,\omega_q$ be corresponding {\em Gauss weights}.  
We denote this quadrature rule, i.e., the set of ordered pairs
$(r_1,\omega_1),\ldots,(r_q,\omega_q)$, by the symbol $\Q$.
Then 
in place of \eref{Kij2} we take
\begin{equation}
K(i,j)=
\sum_{t=1}^m\sum_{\nnu=1}^q\nabla_x \pi_i(\phi_t(r_\nnu))\cdot
\theta(\phi_t(r_\nnu))\nabla_x\pi_j(\phi_t(r_\nnu))
\det(\nabla\phi_t(r_\nnu))\omega_\nnu, \label{Kijq}
\end{equation}
in which $\nabla_x\pi_i(\phi_t(r_\nnu))$ is evaluated by substituting
$z=r_\nnu$ into the right-hand side of \eref{gradpi} and similarly for 
$\nabla_x\pi_j(\phi_t(r_\nnu))$.  Symmetry and sparsity of $K$ follow
for the same reason as for $K^{\rm exact}$; positive semidefiniteness will
follow from results in the next two sections (under some additional
assumptions to be made).

We close this section with a summary of the notation introduced
thus far.  Integers that define the size of the
computation are:
\begin{eqnarray*}
d&=&\mbox{space dimension of \eref{bvp}},\\
p&=&\mbox{polynomial order of the finite element method}, \\
l&=&\mbox{number of reference nodes}=
\left\{\begin{array}{ll}
(p+1)(p+2)/2 & \mbox{if $d=2$}, \\
(p+1)(p+2)(p+3)/6 & \mbox{if $d=3$},
\end{array}
\right. \\
m&=&\mbox{number of elements}, \\
n&=&\mbox{number of non-Dirichlet real-space nodes}, \\
n'&=&\mbox{total number of real-space nodes}, \\
q&=&\mbox{number of Gauss points in the quadrature rule $\Q$}.
\end{eqnarray*}

Sets of points in $\R^d$ include
\begin{eqnarray*}
z_1,\ldots,z_l &=& \mbox{reference nodes}, \\
\zeta_{1,1},\ldots,\zeta_{m,l} &=&\mbox{real-space nodes (local numbering)}, \\
w_1,\ldots,w_{n'} &=&\mbox{real-space nodes (global numbering)}.
\end{eqnarray*}

The quadrature rule $\Q$ is defined by
\begin{eqnarray*}
r_1,\ldots,r_q&=&\mbox{Gauss points}, \\
\omega_1,\ldots,\omega_q&=&\mbox{Gauss weights}. 
\end{eqnarray*}

Important domains are:
\begin{eqnarray*}
T_0 &=&\mbox{the reference element}, \\
\Omega &=& \mbox{the domain of \eref{bvp}}, \\
\phi_1(T_0),\ldots,\phi_m(T_0) &=& \mbox{the elements}, \\
\tilde\Omega&=& \mbox{the approximation to $\Omega$ given by 
$\phi_1(T_0)\cup\cdots\cup\phi_m(T_0)$}.
\end{eqnarray*}

Important functions are:
\begin{eqnarray*}
\phi_1,\ldots,\phi_m &=& \mbox{element mapping functions ($T_0\rightarrow\tilde\Omega$)}, \\
N_1,\ldots,N_l &=& \mbox{shape functions ($T_0\rightarrow \R$)}, \\
\theta &=&\mbox{conductivity ($\tilde\Omega\rightarrow\R$)}, \\
LG &=&\mbox{local-to-global index mapping $\{1,\ldots,m\}\times\{1,\ldots,l\}
\rightarrow\{1,\ldots,n'\}$}, \\
\pi_1,\ldots,\pi_{n'} &=& \mbox{basis functions ($\tilde\Omega\rightarrow\R$)}.
\end{eqnarray*}

Finally, scalar quantities to be introduced in the next section, but which
are included here for completeness, include
\begin{eqnarray*}
m_\Q,M_\Q&=&\mbox{min and max quadrature weights}, \\
\hat \theta(\T)&=&\mbox{intra-element conductivity variation}, \\
\alpha_1,\ldots,\alpha_m&=& \mbox{maximum compression in elements}, \\
\beta_1,\ldots,\beta_m&=& \mbox{maximum stretch in elements}, \\
\kappa_1(\T),\kappa_2(\T)&=&\mbox{mesh quality measures}, \\
\sigma_{\Q,p},\tau_{\Q,p}&=&\mbox{max and min singular values of $S_{\Q,p}$}.
\end{eqnarray*}

This section presented in detail the construction of $K$, the
finite element stiffness matrix, but omitted the construction of $\f$,
the right-hand side of the linear system $K\x=\f$
under consideration.
The construction of $\f$ is not relevant to the main result of
this paper, but for the sake of completeness, we give a very brief overview
and refer the reader to \cite{Johnson} for the details.  
The steps involved in obtaining a formula for the $i$th
entry  $\f$ are as follows.
Multiply
the PDE on both sides by $\pi_i$, integrate over $\Omega$, and
apply integration by parts.  Applying integration by parts causes a boundary
term involving Neumann data $g$ to enter the equation.
Write $u$ as a linear combination of $\pi_j$'s
for $j=1,\ldots,n'$.  The coefficients of this linear combination are
the entries of $\x$ for $j=1,\ldots,n$ and are known from the Dirichlet
boundary condition for $j=n+1,\ldots,n'$.  Finally, gather terms not involving
entries of $\x$ together on the right-hand side and apply quadrature rules
to obtain the $i$th entry of $\f$.
Once the linear system $K\x=\f$ is solved, the approximation of the
unknown function $u$ of \eref{bvp}
is recovered as a linear combination of the $\pi_j$'s using entries of
$\x$ and the Dirichlet conditions as coefficients.

\section{Condition numbers and assumptions}

Our main theorem about matrix approximation, which is Theorem~\ref{factthm}
below, depends on several scalars associated with the finite element
method and the problem at hand
that we define in this section.  In addition, this section
states further assumptions
about the problem and method.

Two constants appearing in our bound
are $m_\Q$ and $M_\Q$, which we define to
be the minimum and maximum weights in the quadrature rule, i.e.,
\begin{equation}
m_\Q=\min(\omega_1,\ldots\omega_q);\; M_\Q=\max(\omega_1,\ldots,\omega_q).
\label{mqdef}
\end{equation}

\begin{assumption}
The quadrature weights are positive, i.e., $m_\Q>0$.
\label{posweight}
\end{assumption}

There is some loss of generality with this assumption because a few
popular finite element quadrature schemes (but certainly not all) use
negative weights \cite{Cools}.

\begin{assumption}
The quadrature scheme is exact for polynomials
of degree up to $2p-2$,
i.e., if $\psi:T_0\rightarrow \R$ is a polynomial of degree $2p-2$ or
less in the $d$ coordinates of $T_0$, then
$$\sum_{\nnu=1}^q \psi(r_\nnu)\omega_\nnu=
\int_{T_0}\psi(z)\,dz.$$
\label{exact2p}
\end{assumption}
This assumption is quite reasonable since it is usually required anyway
for accurate solution by finite element analysis: one wants accurate
quadrature of $\nabla\pi_i\cdot\nabla\pi_j$.

We now define a $dq\times(l-1)$ matrix $S_{\Q,p}$ according to the
following formula:
\begin{equation}
S_{\Q,p} =
\left(
\begin{array}{cccc}
\nabla N_2(r_1) & \nabla N_3(r_1)&\cdots  & \nabla N_l(r_1) \\
\vdots & & & \vdots \\
\nabla N_2(r_q) & \nabla N_3(r_q) & \cdots & \nabla N_l(r_q)
\end{array}
\right).
\label{sqpdef}
\end{equation}
Although we introduce this matrix in this section
in order to define  two associated scalars, 
the motivation for this definition will be postponed until
Section~\ref{reformsec}. 
Before
defining these scalars, we require the following lemma.

\begin{lemma}
Under Assumption~$\ref{exact2p}$, matrix $S_{\Q,p}$ has 
full column rank.
\label{sqplem}
\end{lemma}

\begin{proof}
Let $v\in\R^{l-1}$ be chosen so that $S_{\Q,p}v=\bz$.
Define the function $U:\T_0\rightarrow \R$ according to the formula
$U=v_1N_2+v_2N_3+\cdots+v_{l-1}N_l$.  Observe that, by definition of
$S_{\Q,p}$, the first $d$ entries of $S_{\Q,p}v$ are precisely
$\nabla U(r_1)$, and the next $d$ entries are $\nabla U(r_2)$, etc.  
Therefore, $\nabla U$ vanishes identically at $r_1,\ldots,r_q$.
Let $\psi:T_0\rightarrow\R$ be defined by $\psi=\nabla U \cdot \nabla U$.
Then $\psi$ has degree at most $2p-2$, is a nonnegative function,
and also vanishes identically at $r_1,\ldots,r_q$.  
By Assumption~\ref{exact2p}, this means that $\int_{T_0}\psi=0$.  But
since $\psi$ is nonnegative-valued, we conclude that $\psi\equiv 0$.
Therefore, $\nabla U\equiv 0$ also, and thus $U$ must be a constant
function on $T_0$.  But the definition of $U$ omits the $N_1$ term, and
therefore $U(z_1)=0$.  Since $U$ is constant, this means that $U\equiv 0$
on all of $T_0$, and in particular, $U(z_2)=\cdots=U(z_l)=0$.  But 
$U(z_2)=v_1$, $U(z_3)=v_2$, \ldots, $U(z_l)=v_{l-1}$ by construction
of $U$.  Therefore, the entries of $v$ are all zeros.  Since $v$ was
an arbitrary vector in $N(S_{\Q,p})$, this argument proves that 
this nullspace contains only the $\bz$ vector, hence $S_{\Q,p}$ has
full column rank.
\end{proof}

The next constants that appear in our theorem are $\sigma_{\Q,p}$
and $\tau_{\Q,p}$, defined as follows:
\begin{equation}
\sigma_{\Q,p}=\sigma_{\max}(S_{\Q,p});\; \tau_{\Q,p}=\sigma_{\min}(S_{\Q,p}).
\end{equation}
These constants depend only on $p$ and the quadrature scheme.  
It follows from the lemma that both are positive.

The remaining scalars and assumptions in this section
pertain to the mesh.  Define for $t=1,\ldots,m$,
\begin{eqnarray}
\alpha_t&=&\max\{\Vert \nabla\phi_t(r_1)^{-1}\Vert,\ldots,\Vert \nabla\phi_t(r_q)^{-1}\Vert\},
\label{alphadef} \\
\beta_t&=&\max\{\Vert \nabla \phi_t(r_1)\Vert,\ldots,\Vert \nabla \phi_t(r_q)\Vert\}.
\label{betadef}
\end{eqnarray}
The norms in these equations are matrix 2-norms: $\Vert B\Vert
\equiv\sigma_{\max}(B)$.
Finally, for the whole mesh, a quality measure is
\begin{equation}
\kappa_1(\T)=\max_{t=1,\ldots,m} \alpha_t\beta_t.
\label{qdef}
\end{equation}

Although the $\alpha$'s and $\beta$'s are subscripted only by $t$,
it is clear from the definition that they also depend on $\Q$.  On
the other hand, it is possible to get a $\Q$-independent definition
of these by simply taking the upper bounds similar to \eref{alphadef},
\eref{betadef}
over all $z\in T_0$ instead of just $r_1,\ldots,r_q$.  
Unfortunately, for higher order elements ($p>1$), there is no
simple method to compute $\max_{z\in T_0}\Vert \nabla\phi_t(z)^{-1}\Vert$; 
a technique for obtaining an upper bound appears in \cite{Vava:bezier}.

It should be noted that $\kappa_1(\T)\ge 1$ since for any $z$,
$\Vert \nabla\phi_t(z)^{-1}\Vert\cdot\Vert \nabla \phi_t(z)\Vert\ge 1$.
If all the elements are well-shaped, i.e., not too distorted
when compared to the reference element, then $\kappa_1(\T)$
will not be much larger than 1.

The second mesh quality measure is
\begin{equation}
\kappa_2(\T)=\max_{t=1,\ldots,m} \frac
{\max_{\nnu=1,\ldots,q} \det(\nabla \phi_t(r_\nnu))}
{\min_{\nnu=1,\ldots,q} \det(\nabla \phi_t(r_\nnu))}.
\label{qdef2}
\end{equation}
This quantity measures the maximum over elements of
the variation in volumetric distortion over the element.  This
may be regarded as a measure of how much elements depart from
linearity (flatness).
Measure $\kappa_2$ is not completely independent from $\kappa_1$
as the following argument shows.
It follows from Hadamard's inequality that
$$\alpha_t^{-d} \le \det(\nabla \phi_t(r_\nnu)) \le \beta_t^d$$
and therefore
$$\frac
{\max_{\nnu=1,\ldots,q} \det(\nabla \phi_t(r_\nnu))}
{\min_{\nnu=1,\ldots,q} \det(\nabla \phi_t(r_\nnu))}
\le (\alpha_t\beta_t)^d,$$
hence $\kappa_2(\T)\le \kappa_1(\T)^d$.  This bound is not likely
to be tight in practice (see our computation results in 
Table~\ref{condtable}), 
so we prefer to distinguish the roles
of $\kappa_1$ and $\kappa_2$.

The final two assumptions are qualitative in nature and are meant to
indicate cases in which our method is expected to work well. Violation
of these assumptions does not invalidate our theorem but merely indicates
that our bound on $\kappa(K,\bar K)$ may grow large.

\begin{assumption}
Mesh quality measures $\kappa_1(\T),\kappa_2(\T)$ are
not too large.
\label{kappaassum}
\end{assumption}

In the case of linear elements ($p=1$), 
$\kappa_1(\T)$ is proportional to the reciprocal of the minimum
angle (in two dimensions) or solid angle (in three dimensions) of the mesh
and $\kappa_2(\T)\equiv 1$.

Assumption~\ref{kappaassum} does not imply that
the elements are of a uniform size: a uniform rescaling of element
$t$ does not affect the product $\alpha_t\beta_t$.

The final constant and assumption pertain to the conductivity
field.  Define
\begin{equation}
\hat\theta(\T) = \max_{t=1,\ldots,m} \frac
{\max_{\nnu=1,\ldots,q}\theta(\phi_t(r_\nnu))}
{\min_{\nnu=1,\ldots,q}\theta(\phi_t(r_\nnu))}.
\label{hathetadef}
\end{equation}
In other words, $\hat\theta(\T)$ measures the maximum intra-element variation
of the conductivity field $\theta$.

\begin{assumption}
Intra-element conductivity variation $\hat\theta(\T)$ is not too large.
\end{assumption}

This assumption implies that if there are huge jumps in $\theta$ in
the domain (e.g., because one is modeling a domain composed of two
materials with vastly different conductivities), then the mesh
boundaries should be aligned with the conductivity jumps.  
If the mesh boundaries are aligned with the conductivity jumps, then
no element will have large variation in $\theta$ among its Gauss points
and thus $\hat\theta(\T)$ will not be large.

Note that, as in the case of $\kappa_1(\T)$ and $\kappa_2(\T)$, 
scalar $\hat\theta(\T)$ also
depends on the quadrature rule, but this dependence can be eliminated
by overestimating $\hat\theta(\T)$ as 
$$\max_{t=1,\ldots,m} \frac
{\max_{z\in T_0}\theta(\phi_t(z))}
{\min_{z\in T_0}\theta(\phi_t(z))}.$$

\section{The matrix approximation}
\label{sec:matrixapprox}

Our main matrix factorization result is summarized by
the following theorem whose proof is 
explained in upcoming sections.

\begin{theorem}
\label{factthm}
Let $K$ be defined by \eref{Kijq} above, and
let Assumptions~$\ref{postheta}$--$\ref{exact2p}$ hold.
Then $K$ may be factored as $A^TJ^TDJA$, where 
\begin{itemize}
\item
$A$ is an $(l-1)m\times n$
reduced node-arc incidence matrix of a certain multigraph, 
\item
$J$ is a $dqm\times(l-1)m$ matrix
that is well conditioned in the sense that
\begin{equation}
\sigma_{\max}(J)\le \sigma_{\Q,p},
\label{Jboundmax}
\end{equation}
and
\begin{equation}
\sigma_{\min}(J)\ge \tau_{\Q,p}/\kappa_1(\T).
\label{Jboundmin}
\end{equation}
\item
$D$ is a $dqm\times dqm$ positive definite diagonal matrix.
\end{itemize}
Further, $J^TDJ$ may be refactored as $\bar{D}^{1/2}H\bar{D}^{1/2}$
where $\bar D$ is a $(l-1)m\times (l-1)m$ SPD diagonal
matrix and $H$ is a $(l-1)m\times (l-1)m$ SPD matrix
whose condition number is bounded as follows:
\begin{equation}
\kappa(H) \le \hat\theta(\T)\kappa_1(\T)^2\kappa_2(\T)\cdot 
\frac{M_\Q\sigma_{\Q,p}^2}{m_\Q\tau_{\Q,p}^2}.
\label{hbd}
\end{equation}
\end{theorem}

This theorem can now be combined with Theorem~\ref{thm:condno}
to obtain a good approximation $\bar K$ to the stiffness
matrix $K$.  In particular, we take $\bar K=A^T\bar{D}A$.  Note
that this matrix, being a weighted graph laplacian, is symmetric
and diagonally dominant.
Then in the context of Theorem~\ref{thm:condno}, $V=A^T\bar{D}^{1/2}$.  
The theorem thus implies that the condition
number of $\bar K$ with respect to $K$ depends only on the condition
number of $H$ for which we have a good bound \eref{hbd}.  Then $\bar K$
can be preconditioned using techniques in the previous literature.
The complete description of the algorithm
is given below in Section~\ref{algo}.

The approximation bound clearly depends on the quadrature rule,
mesh quality, and intra-element variation of $\theta$.  
It is also informative to make a list of quantities on which
$\kappa(\bar K, K)$ does {\em not} depend:
\begin{itemize}
\item
The number of nodes or elements,
\item
The size of the elements,
\item
Variation in the size of elements (i.e., gradation of the mesh),
\item
The shape of $\Omega$,
\item
The conductivity field $\theta$
(provided that the mesh respects internal boundaries
where large conductivity jumps occur).
\end{itemize}
Although the shape of $\Omega$ does not directly affect 
$\kappa(\bar K,K)$, naturally, it affects the mesh generation procedure
which in turn affects mesh quality.  A particular case that requires
further consideration is when $\Omega$ 
has a sharp corner; a sharp corner of $\Omega$ forces
the mesh to contain at least one
element with a sharp corner, thus causing
$\kappa_1(\T)$ to be large.  Even in this case, however, it may
be possible to obtain an upper bound on $\kappa(\bar K, K)$ independent
of this sharp angle using the analysis of
Phillips and Miller (see Section~\ref{sec:notesaddedinrevision}).

\section{The first factorization}
\label{reformsec}
In this section we develop the first factorization of $K$
stated in Theorem~\ref{factthm}. 
This factorization is related to one proposed in \cite{Vava:fe}.
We
state the main result of this section as a lemma.

\begin{lemma}
Matrix $K$ defined by \eref{Kijq} can be factored as
$K=A^TS^TR^TDRSA$ where the factors are as follow.
\begin{itemize}
\item
Define matrix $A\in \R^{(l-1)m\times n}$ to be a sparse matrix
all of whose entries are $-1$, $0$ or $1$.  In more
detail, $A$ is written in block form
$$
A=\left(
\begin{array}{c}
A_1 \\
\vdots \\
A_m
\end{array}
\right),
$$
where $A_t$ is $(l-1)\times n$ for each element $t=1,\ldots,m$.
The columns are indexed $1,\ldots, n$ in correspondence
with nodes $w_1,\ldots,w_n$.  Row $\mu-1$ of $A_t$ has a `$1$' in 
column $LG(t,\mu)$  for $\mu=2,\ldots,l$ and a `$-1$' in column 
$LG(t,1)$.  
Thus, most rows of $A$ have exactly two nonzero entries.
If $LG(t,\mu)>n$ (i.e., node $\zeta_{t,\mu}$ lies in 
the Dirichlet boundary $\tilde\Gamma_1$),
then the `$1$' entry is omitted.  Similarly, if $LG(t,1)>n$, then
the `$-1$' entry is omitted.  For this reason, a few rows of $A$ have just one
nonzero entry or none at all.
\item
Define $S\in \R^{dqm\times (l-1)m}$ by
\begin{equation}
\left.
S=\left(
\begin{array}{ccc}
S_{\Q,p}  & & \\
 & \ddots &\\
& & S_{\Q,p}
\end{array}
\right)
\right\} \mbox{$m$ times,}
\label{Sdef}
\end{equation}
where $S_{\Q,p}$ was defined by \eref{sqpdef}.
\item
Define block diagonal  $R\in \R^{dqm\times dqm}$  by
\begin{equation}
R=\left(
\begin{array}{ccc}
R_1 & & \\
& \ddots & \\
& & R_m
\end{array}
\right),
\label{Sdef2}
\end{equation}
where each of $R_1,\ldots,R_m\in \R^{qd\times qd}$ 
is itself block diagonal and
given by
\begin{equation}
R_t=\alpha_t^{-1}\left(
\begin{array}{ccc}
\nabla\phi_t(r_1)^{-T} & & \\
& \ddots & \\
& & \nabla\phi_t(r_q)^{-T}
\end{array}
\right).
\label{rtdef}
\end{equation}
The scalar $\alpha_t$ used here was defined by \eref{alphadef}.
\item
Finally, $D \in \R^{qdm\times qdm}$ is positive definite diagonal and
given by
\begin{equation}
D=\left(
\begin{array}{ccc}
D_1 & & \\
& \ddots & \\
& & D_m
\end{array}
\right),
\label{Ddef}
\end{equation}
where $D_t\in \R^{qd\times qd}$, $t=1,\ldots,m$, is given by
\begin{equation}
D_t=\alpha_t^2 \left(
\begin{array}{ccc}
\theta(\phi_t(r_1))\det(\nabla\phi_t(r_1))
\omega_1I & & \\
& \ddots & \\
& & \theta(\phi_t(r_q))\det(\nabla\phi_t(r_q))
\omega_qI  
\end{array}
\right)
\label{dmatdef}
\end{equation}
in which $I$ denotes the $d\times d$ identity matrix.
\end{itemize}
\label{firstfaclem}
\end{lemma}

\noindent
{\em Remarks.}
Intuitively, this factorization of $K$ decomposes finite element
analysis into natural ingredients:
$A$ encodes the combinatorial connectivity
of the mesh, $R$ encodes the geometry of the mesh, $S$ encodes the
quadrature points, and $D$ encodes the
quadrature weights and conductivity.  The scaling factor $\alpha_t$,
which is necessary for our analysis, cancels out between $D_t$ and $R_t$.

Note that the positive definiteness of $D$ follows from 
Assumptions~\ref{postheta}, \ref{posdetphi}, and \ref{posweight}.
Note also that $A$ is a reduced
node-arc incidence matrix of a multigraph defined on the nodes
of $\T$.  ({\em Multi}graph indicates that more than one edge may connect
a particular pair of vertices.)
Each element gives rise to $l-1$ arcs in the graph.  
In particular, for each element $t=1,\ldots,m$, there is an arc
joining each of its nodes $2,\ldots,l$ to node 1.
Columns corresponding to nodes 
of $\tilde \Gamma_1$ are omitted.
(This is what is meant by ``reduced.'')

\begin{proof}
Let $\bxi$ be an arbitary vector in $\R^n$, and let 
$u:\tilde\Omega\rightarrow\R$
be defined by
$$u=\sum_{i=1}^n \xi_i\pi_i.$$
For an element $t$, define function $U_t:T_0\rightarrow\R$ by
\begin{equation}
U_t=\sum_{\mu=1}^l \xi_{LG(t,\mu)}N_\mu.
\label{UTdef}
\end{equation}
With these definitions, it is clear that
$u\circ\phi_t=U_t$.  In \eref{UTdef}, we follow the convention that
$\xi_i\equiv 0$ in the case that $i\in\{n+1,\ldots,n'\}$.

The assembled stiffness matrix is defined by \eref{Kijq} so that
\begin{eqnarray}
\bxi^TK\bxi & = & \sum_{i=1}^n\sum_{j=1}^n\xi_i\xi_jK(i,j) \nonumber \\
&=&
\sum_{i=1}^n\sum_{j=1}^n\xi_i\xi_j
\sum_{t=1}^m\sum_{\nnu=1}^q\nabla_x \pi_i(\phi_t(r_\nnu))\cdot
\theta(\phi_t(r_\nnu))\nabla_x\pi_j(\phi_t(r_\nnu))
\det(\nabla \phi_t(r_\nnu))\omega_\nnu \nonumber
\\
&=&
\sum_{t=1}^m\sum_{\nnu=1}^q
\sum_{i=1}^n\xi_i \nabla_x \pi_i(\phi_t(r_\nnu))\cdot
\theta(\phi_t(r_\nnu))
\sum_{j=1}^n\xi_j \nabla_x\pi_j(\phi_t(r_\nnu))
\det(\nabla\phi_t(r_\nnu))\omega_\nnu \nonumber
\\
&=&
\sum_{t=1}^m\sum_{\nnu=1}^q
\nabla_x u(\phi_t(r_\nnu))\cdot
\theta(\phi_t(r_\nnu))
\nabla_x u(\phi_t(r_\nnu))
\det(\nabla\phi_t(r_\nnu))\omega_\nnu \nonumber
\\
&=&
\sum_{t=1}^m\sum_{\nnu=1}^q
\nabla_x U_t(r_\nnu)\cdot
\theta(\phi_t(r_\nnu))
\nabla_x U_t(r_\nnu)
\det(\nabla\phi_t(r_\nnu))\omega_\nnu \nonumber
\\
&=&
\sum_{t=1}^m\sum_{i=1}^q (\nabla\phi_t(r_\nnu)^{-T}\nabla_z U_t(r_\nnu))
\cdot\theta(\phi_t(r_\nnu)) \nonumber \\
& & \qquad\qquad
(\nabla\phi_t(r_\nnu)^{-T}\nabla_z U_t(r_\nnu))
\det(\nabla\phi_t(r_\nnu))\omega_\nnu \label{vDv1}\\
&=&
\v^TD\v. \label{vDv2}
\end{eqnarray}
In \eref{vDv2}, we
have used the matrix $D$ defined by \eref{Ddef}.  We
have also introduced the vector $\v\in\R^{dqm}$ defined
in block fashion as follows.  Write $\v=[v_1;\cdots;v_m]$
where $v_t\in\R^{dq}$ is itself composed of blocks
$v_t=[v_{t,1};\cdots;v_{t,q}]$.  Here $v_{t,\nnu}\in\R^d$
is defined to be 
\begin{equation}
v_{t,\nnu}=(\nabla\phi_t(r_\nnu))^{-T}\nabla_z U_t(r_\nnu)/\alpha_t.
\label{vtnu}
\end{equation}
It is clear by construction of $D$ and $\v$ that $\v^TD\v$
is equal to the expression in \eref{vDv1}.

Next, we claim that 
$\v=RSA\bxi$.  Focus on the block corresponding to a
particular element $t\in\{1,\ldots,m\}$ in which case we
must show that $v_t=R_tS_{\Q,p}A_t\bxi$.
The product $A_t\bxi$ yields a vector with $(l-1)$ entries
that
contains  finite differences
of entries of $\bxi$.  Specifically, the $\mu-1$ entry is
$\xi_{LG(t,\mu)}-\xi_{LG(t,1)}$ for $\mu=2,\ldots,l$.

By definition of $S_{\Q,p}$ in \eref{sqpdef}, 
it follows that the $\nnu$-block
of entries ($\nnu=1,\ldots,q$) of 
$S_{\Q,p}A_t\bxi$ is the $d$-vector
$\nabla_z \bar U_t(r_\nnu)$, where 
$$\bar U_t=(\xi_{LG(t,2)}-\xi_{LG(t,1)})N_2+\cdots +
(\xi_{LG(t,l)}-\xi_{LG(t,1)})N_l.$$
Comparing this equation to \eref{UTdef} indicates that
$\bar U_t$ and $U_t$ differ by 
$$\xi_{LG(t,1)}(N_1+\cdots+N_l).$$
This latter quantity, however, is a constant function
(because $N_1+\cdots+N_l$ is identically 1, a property that
follows from \eref{deltafunc} and the fact that 1 is a polynomial
of degree at most $p$ and hence must be expressable as a sum of
$N_\mu$'s).  Therefore, $U_t$ and $\bar U_t$ have the same gradient.
We conclude that the $\nnu$ block of $S_{\Q,p}A_t\bxi$ must equal
$\nabla_z U_t(r_\nnu)$.

Combining this equation with \eref{rtdef} shows that 
$$R_tS_{\Q,p}A_t\bxi=(\nabla\phi_t(r_\nnu))^{-T}\nabla_z 
U_t(r_\nnu)/\alpha_t,$$
and hence is equal to $v_{t,\nnu}$ as defined by \eref{vtnu}.
This concludes the proof that $\v=RSA\bxi$.

Since $\v^TD\v=\bxi^TK\bxi$, this means
$\bxi^TA^TS^TR^TDRSA\bxi=\bxi^TK\bxi$.  Since this holds for
all $\bxi\in\R^n$, and since a symmetric matrix $C$ is uniquely determined
by the mapping $\bxi\mapsto\bxi^TC\bxi$, this means that
$K=A^TS^TR^TDRSA$, which concludes the proof of the lemma.
\end{proof}

Now we analyze the singular values of $J$
to finish the first part of Theorem~\ref{factthm}.
Let us define $J=RS$ so that $K=A^TJ^TDJA$ as claimed.
The block structures of $R$ and $S$ induce a corresponding
block structure on $J$:
$$J=\left(
\begin{array}{ccc}
J_1 & & \\
& \ddots & \\
& & J_m
\end{array}
\right)
$$
where 
\begin{equation}
J_t=R_tS_{\Q,p}
\label{Jtdef}
\end{equation}
for all $t=1,\ldots,m$.
Because of this structure, the maximum singular value of $J$ is the maximum
singular value
among any of its blocks and similarly for its minimum singular value.
Since in general $\sigma_{\max}(AB)\le \sigma_{\max}(A)\sigma_{\max}(B)$,
\begin{eqnarray}
\sigma_{\max}(J_t) & \le& \sigma_{\max}(R_t)\sigma_{\max}(S_{\Q,p}) 
\nonumber \\
&=& \sigma_{\max}(\diag(\nabla\phi_t(r_1)^{-T},\ldots, \nabla\phi_t(r_q)^{-T}))
\sigma_{\Q,p}/\alpha_t \nonumber \\
&=&\max_{\nnu=1,\ldots,q}\sigma_{\max}(\nabla\phi_t(r_\nnu)^{-1})
\sigma_{\Q,p}/\alpha_t \nonumber\\
& =  & \sigma_{\Q,p}. \label{Jtsigmamax}
\end{eqnarray}
The last line follows from \eref{alphadef} and
establishes \eref{Jboundmax}.

Since $\sigma_{\min}(AB)\ge \sigma_{\min}(A)\sigma_{\min}(B)$ for two
matrices $A,B$ with full column rank,
\begin{eqnarray}
\sigma_{\min}(J_t)&\ge& \sigma_{\min}(R_t)\sigma_{\min}(S_{\Q,p}) 
\nonumber\\
&=&
\sigma_{\min}(\diag(\nabla\phi_t(r_1)^{-T},\ldots, \nabla\phi_t(r_q)^{-T}))
\tau_{\Q,p}/\alpha_t \nonumber \\
&= &\min_{\nnu=1,\ldots,q}\sigma_{\min}(\nabla\phi_t(r_\nnu)^{-1})
\tau_{\Q,p}/\alpha_t \nonumber \\
&=& \min_{\nnu=1,\ldots,q}(1/\sigma_{\max}(\nabla\phi_t(r_\nnu)))
\tau_{\Q,p}/\alpha_t\nonumber \\
&=& (1/\beta_t)\cdot\tau_{\Q,p}/\alpha_t \nonumber\\
&\ge& \tau_{\Q,p}/\kappa_1(\T). \label{Jtsigmamin}
\end{eqnarray} 
For the last line, we used the fact that $\alpha_t\beta_t\le \kappa_1(\T)$,
which follows from \eref{qdef}.
This establishes \eref{Jboundmin} and
concludes the proof of the first factorization in 
Theorem~\ref{factthm}.

Our factorization $K=A^TJ^TDJA$ is reminiscent of
one proposed by Argyris
\cite{Argyris} 
of the form $K=\tilde AP \tilde A^T$,
which he calls
the ``natural factorization.''
In Argyris's factorization, however, the matrix $\tilde A$ has all $+1$
and $0$ entries and therefore is not a node-arc incidence matrix.
The purpose of Argyris's matrix $\tilde A$ is to assemble the element
stiffness matrices, which constitute the blocks of the block-diagonal
matrix $P$.

\section{The second factorization}
\label{sec:secondfac}

In this section we prove the second part of Theorem~\ref{factthm}.
We state the existence of the second factorization in the form of a lemma.
\begin{lemma}
Let $J$ and $D$ be defined as in Lemma~$\ref{firstfaclem}$.  
Then $J^TDJ$ can be refactored as $\bar{D}^{1/2}\bar{J}^T\bar{J}\bar{D}^{1/2}$, where
\begin{itemize}
\item
Matrix $\bar{D}\in\R^{m(l-1)\times m(l-1)}$ is positive definite diagonal
and is written in block form
\begin{equation}
\bar{D}=m_\Q\left(
\begin{array}{ccc}
f_1g_1\alpha_1^2I& & \\
& \ddots & \\
& & f_mg_m\alpha_m^2I
\end{array}
\right)
\label{barddef}
\end{equation}
where $I$ is the $(l-1)\times(l-1)$ identity matrix.  In this
formula, $f_t$ is the minimum value of $\theta$ over
Gauss points of element $t$:
\begin{equation}
f_t=\min_{\nnu=1,\ldots,q} \theta(\phi_t(r_\nnu)),
\label{ftdef}
\end{equation}
and similarly, $g_t$ is the minimum value of $\det(\nabla \phi_t)$:
\begin{equation}
g_t =\min_{\nnu=1,\ldots,q} \det(\nabla \phi_t(r_\nnu)).
\label{gtdef}
\end{equation}
\item
Matrix $\bar J\in \R^{dqm\times(l-1)m}$ is defined by
\begin{equation}
\bar J =D^{1/2}J\bar{D}^{-1/2}.
\label{barjdef}
\end{equation}
\end{itemize}
\end{lemma}

\begin{proof}
The fact that $\bar{D}^{1/2}{\bar J}^T\bar{J}\bar{D}^{1/2}=J^TDJ$
follows as an immediate consequence of \eref{barjdef} regardless
of how we have defined $\bar{D}$.
\end{proof}

Next we analyze the condition number of $H$ to finish proving the
second part of Theorem~\ref{factthm}.
We define $H={\bar J}^T\bar J$
so that $J^TDJ=\bar{D}^{1/2}H\bar{D}^{1/2}$.  We
must estimate
the singular values of $\bar J$, which are the square roots of the
eigenvalues of $H$.

Let the diagonal block of $\bar J$ associated with 
element $t$ be denoted $\bar J_t$ for $t=1,\ldots,m$.
Because of the block structure, the maximum and minimum singular values
for $\bar J$ are the maximum and minimum singular values among the
blocks $\bar J_t$ which may be written
\begin{eqnarray*}
\bar J_t&=& D_t^{1/2}J_t
f_t^{-1/2}g_t^{-1/2}m_\Q^{-1/2}\alpha_t^{-1} \\
&=& \hat D_t J_t
\end{eqnarray*}
where $$\hat D_t = f_t^{-1/2}g_t^{-1/2}m_\Q^{-1/2}
\alpha_t^{-1}D_t^{1/2}.$$
Examining the constituent parts of $D_t$ given by
\eref{dmatdef}, observing that the $\alpha_t^2$ in \eref{dmatdef}
is cancelled by the $\alpha_t^{-1}$ in the preceding equation
and applying the inequalities 
$$1\le \omega_\nnu/m_\Q\le M_\Q/m_\Q$$
for $\nnu=1,\ldots,q$ (see \eref{mqdef}),
$$1\le \theta(\phi_t(r_\nnu)) f_t^{-1} \le
\frac{\max_{\nnu=1,\ldots,q} \theta(\phi_t(r_\nnu))}{\min_{\nnu=1,\ldots,q}\theta(\phi_t(r_\nnu))}$$
(see \eref{ftdef}), and 
$$1\le \det(\nabla\phi_t(r_\nnu))g_t^{-1}\le 
\frac{\max_{\nnu=1,\ldots,q}\det(\nabla \phi_t(r_\nnu))}{\min_{\nnu=1,\ldots,q}\det(\nabla\phi_t(r_\nnu))}$$
(see \eref{gtdef}), we conclude that
the diagonal entries of $\hat D_t$ satisfy
$$1\le \hat D_t(i,i)  \le
\left(\frac
{\max_{\nnu}\theta(\phi_t(r_\nnu))}
{\min_{\nnu}\theta(\phi_t(r_\nnu))}
\right)^{1/2}
\left(\frac
{\max_{\nnu}\det(\nabla\phi_t(r_\nnu))}
{\min_{\nnu}\det(\nabla\phi_t(r_\nnu))}
\right)^{1/2}\left(\frac{M_\Q}{m_\Q}\right)^{1/2}.$$
By \eref{qdef2} and \eref{hathetadef}, the first two quantities on the
right-hand side of the preceding equation are bounded by
$\hat\theta(\T)^{1/2}$ and $\kappa_2(\T)^{1/2}$ respectively.
Thus, it is apparent that for each $i$,
$$1\le D_t(i,i)\le  (\hat\theta(\T)\kappa_2(\T)M_\Q/m_\Q)^{1/2}.$$
Since $\bar{J}_t=\hat{D}_tJ_t$, we
can combine the inequalities in the previous line with \eref{Jtsigmamax}
and \eref{Jtsigmamin} to obtain:
\begin{equation}
\sigma_{\max}(\bar J_t) \le \hat\theta(\T)^{1/2}\kappa_2(\T)^{1/2} M_\Q^{1/2}\sigma_{\Q,p}/m_\Q^{1/2},
\label{Jboundmax2}
\end{equation}
and
\begin{equation}
\sigma_{\min}(\bar J_t) \ge \tau_{\Q,p}/\kappa_1(\T).
\label{Jboundmin2}
\end{equation}
\label{Jbound2}

Since $H=\bar{J}^T\bar{J}$ and $\bar J$ has full column rank, 
\eref{Jboundmax2} and \eref{Jboundmin2} prove \eref{hbd}, which concludes
the proof of Theorem~\ref{factthm}.

\section{Preconditioning Strategy Summary}
\label{algo}

The main result of this paper is that the stiffness matrix $K$
of \eref{bvp}
can be approximated by a symmetric diagonally dominant
matrix $\bar{K}$.  In this section we discuss several
approaches for using this approximation to efficiently
solve the linear system $K\x=\f$ by iteration.

The first step in using the approximation to construct $\bar{K}$.
Recall that $\bar{K}=A^T\bar{D}A$, and thus $A$ and $\bar D$ must
be formed.
Forming $A$ means construction of the multigraph
consisting of $l-1$ edges (a star-tree) per element of $\T$.  Matrix
$A$ is then the reduced node-arc incidence matrix of this multigraph
as specified by Lemma~\ref{firstfaclem}.  Construction of
$\bar{D}$ is given by \eref{barddef}.
Notice that computing
$A$ and $\bar{D}$ requires information about the mesh and original
boundary value problem.  In other words, our method is not applicable
(at least not in a straightforward manner)
if the only information about the original problem is the
stiffness matrix $K$.  

Once $\bar{K}$ is on hand, there are several ways to proceed.
The most straightforward is
to  suppose that one has an efficient preconditioned
conjugate gradient solver for systems of the form $\bar{K}\x=\f$.
Efficient preconditioners for
symmetric diagonally
dominant systems were proposed and analyzed in a graph-theoretic
framework by 
Vaidya.  Vaidya's work is described and extended by
\cite{BGHNT,CT03}, and \cite{Gre96} contains a related analysis.
One of the most recent improvements is due
to Spielman and Teng \cite{ST04}, who
propose a graph-based preconditioner whose running time is
$O(n^{5/4})$.  (Spielman and Teng improve this bound to $O(n^{1+\epsilon})$
for any $\epsilon>0$, but the
algorithm corresponding to the improved bound is no longer preconditioned
conjugate gradients.)  Other  techniques proposed in the literature
for symmetric diagonally dominant matrices such as
algebraic multigrid could also be used.

If one has a good preconditioner $M$ for $\bar{K}$, then one
could also use $M$ as a preconditioner directly for $K$.
This is because of the ``triangle inequality'' (see
\cite{BH04}), which states
\begin{equation}
  \sigma(K,M) \le \sigma(K,\bar{K}) \sigma(\bar{K},M).
\label{triangleinequality}
\end{equation}
This shows that the the overall support number 
is bounded by the product of the support numbers in each step of the
approximation $K \approx \bar{K} \approx M$.  In particular.  
$\kappa(M,K)\le \kappa(H)\kappa(M,\bar K)$.  Since the number
of iterations of preconditioned conjugate gradients is bounded by
the square root of the condition number, our analysis shows that
the increase in the number of iterations for solving $K\x=\f$
(compared to solving $\bar K\x=\f$) is at most a factor of
$\sqrt{\kappa(H)}$,
which in turn is bounded by
$$\hat\theta(\T)^{1/2}\kappa_1(\T)\kappa_2(\T)^{1/2}
M_\Q\sigma_{\Q,p}/(m_\Q\tau_{\Q,p}).$$

The asymptotically fastest known iterative algorithm for solving
symmetric diagonally dominant matrices is due to Elkin et al.~\cite{EEST}.
This algorithm, which we denote EEST,
extends Spielman and Teng \cite{ST04} and
requires $O(n_e(\log n)^s)$ time to solve any
diagonally dominant $\bar{K}\x=\f$,
where $n_e$ is the number of nonzero entries in
$\bar K$ and $s$ is some constant.  
Our preconditioner $\bar{K}$  is the Laplacian of a graph with
$ml$ edges, where $m$ is the number of
elements and $l$ is the number of nodes in the reference element
(refer to Lemma~\ref{firstfaclem}), hence $n_e=ml$.
We may assume that $m=O(n)$ for the following
reasons.  First, this relationship holds for all 
automatic finite element mesh generators of which we are aware.
Second, a
violation of the relationship $m=O(n)$ implies that the mesh has
elements of unbounded aspect ratio which in turn means $\kappa_1(\T)$
tends to infinity so that our results
are not as useful anyway. 
Thus, it is safe to assume $n_e\le \mbox{const}\cdot nl$, and we are
regarding $l$ (which depends on $p$, the polynomial degree, and $d$,
the space dimension) as a constant.

The EEST iteration is more complicated than preconditioned conjugate gradients,
and hence it is unclear whether it is possible to simply replace the
diagonally
dominant coefficient matrix 
$\bar K$ by our stiffness matrix $K$ and expect the algorithm
to converge in $O(n(\log n)^s\kappa(H)^{1/2})$ time.  

On the other
hand, one could obtain this running time
$O(n(\log n)^s\kappa(H)^{1/2})$ by using a nested iteration:
the outer iteration is preconditioned conjugate gradients in which $\bar K$
preconditions $K$.  Thus, $O(\kappa(H)^{1/2})$ 
outer iterations are required.
The inner loop is the EEST algorithm to apply the preconditioner,
i.e., to solve $\bar K\x=\f$ in $O(n(\log n)^s)$ iterations.
The EEST algorithm is fairly complex, and
using it in a two-level manner would raise a number of difficulties
associated with termination tests, so it is unlikely to be
practical currently.

This bound of $O(n(\log n)^s\kappa(H)^{1/2})$ on the number of operations is
independent of the condition number of the underlying system because, as
noted above $\kappa(H)$ depends on factors associated with the finite element
method and  on mesh quality measures but not on the conditioning of $K$.

Our approximation scheme can be rewritten on an element-by-element
basis as follows.
For some element index $t\in\{1,\ldots, m\}$,
consider its element stiffness
matrix $K_t = A_t^T J_t^T D_t J_t A_t$, where 
$A_t$ is defined in Lemma~\ref{firstfaclem},
$J_t$ is defined by \eref{Jtdef},
and $D_t$ is defined by \eref{dmatdef}.
The above proof shows that
$\bar{D}_t$, which is the $t$th block
of \eref{barddef},
is a good approximation to $J_t^TD_tJ_t$.
Thus, we let $\bar{K}_t = 
A_t^T \bar{D}_t A_t$.
The overall approximation is $\bar{K} = \sum_{t=1}^m \bar{K}_i$.
Note that the so-called 
splitting lemma frequently used in support-graph
theory (see, e.g., \cite{BH04}) shows that
$$\sigma(K,\bar{K})\le \max(\sigma(K_1,\bar K_1),\ldots,\sigma(K_m,\bar K_m)),$$
so it suffices to analyze the quality of approximation of $\bar{K}_t$
to $K_t$.

Alternatively, we can take a global view as in the above
analysis and write
$\bar{K} = A^T \bar{D} A$, where $A=\left(A_1;A_2; \ldots; A_m \right)$
and $\bar{D} = \diag(\bar{D}_1, \bar{D}_2, \ldots, \bar{D}_m)$. 
To simplify notation, we have adopted the global view
in this paper, but the reader should 
keep in mind that our approximation can take place element by element,
which may be important in an implementation.

\section{Notes added in revision}
\label{sec:notesaddedinrevision}

Since the first version of this paper, there have been three 
related developments of interest.  First, R.~Gupta's Master's thesis
\cite{gupta} gave a geometric method to approximate element stiffness matrices
by diagonally dominant matrices in the $p=1$, $d=2$ case.
In follow up work, Wang and Sarin \cite{WangSarin} 
showed how to parallelize the algorithm.
Computational experiments are presented.  For this particular case
$p=1,d=2$, their result is essentially equivalent to ours.

Second, Phillips and Miller
\cite{PhilMil} have shown that in the $p=1$, $d=2$ case, $\kappa(H)$ is
still small even if some elements have very small angles, provided
that no element has large angles (close to $\pi$).  In our preceding analysis,
the right-hand side of \eref{hbd} would grow large in the presence
of small angles since the factor $\kappa_1(\T)$ would be large.
Wang and Sarin make a similar observation.

In a third development, Avron et al.\ \cite{Avron} have shown
how to get an approximation to element stiffness matrices that is
optimal up to a constant factor.  In other words, for an element
stiffness matrix $K_t$ they find a diagonally dominant matrix $\tilde K_t$
such $\kappa(K_t,\tilde K_t)\le c_1\kappa(K_t, K_t')$, where
$K_t'$ is any other symmetric diagonally dominant matrix of the
correct size and $c_1$ depends only
on $p,d$.
Thus, their preconditioner could  lead to a faster algorithm
than ours since ours is not optimal in this sense.   
Their theory, however, does not subsume ours since
they do not obtain any new bounds on $\kappa(K,\tilde{K})$ comparable
to \eref{hbd}.  Their results are competitive and in some cases better
than an algebraic multgrid method.

\section{Computational tests}
\label{sec:comptest}

There are at least two possible ways to test the effectiveness of
our result: calculate the quality of matrix approximation or measure the
speed of an iterative method.  Since our theory is primarily about the
former issue, and since Avron et al.\ conduct extensive testing on the
latter question, we focus on the first question.

We try out four specific commonly occurring cases: $(d,p)=(2,1),
(2,2),(3,1),(3,2)$, i.e., linear  and quadratic triangles and linear and
quadratic
tetrahedra. 
Six test meshes were tried, which we denote
$\TT_1, \ldots, \TT_6$; the first two are two-dimensional and the
last four are three-dimensional.
For quadrature in the $p=1$ cases, we use the midpoint rule.
For the $p=2$ cases, we use the
symmetric $d+1$-point rule, which is accurate for polynomials up to
degree 2.  Table~\ref{quadtable} gives the quadrature
rules used.

\begin{table}
\caption{Quadrature rules used for the six test cases.
The rules are written in the form $(r_1,\omega_1),\ldots,(r_q,\omega_q)$,
where $r_\nnu$ is a Gauss point in the reference simplex and $\omega_\nnu$
is its weight.  In this table $\xi_1=(10-\sqrt{20})/40$ and $\xi_2=1-3\xi_1$.}
\label{quadtable}
\begin{center}
\begin{tabular}{ll}
\hline
Test & $\Q$ \\
\hline
$\TT_1$ & $((1/3,1/3),1/2)$ \\
$\TT_2$ & $((1/6,1/6),1/6), ((1/6,2/3),1/6), ((2/3,1/6),1/6)$ \\
$\TT_3,\TT_5$ & $((1/4,1/4,1/4),1/6)$ \\
$\TT_4,\TT_6$ & $((\xi_1,\xi_1,\xi_1),1/24),((\xi_1,\xi_1,\xi_2),1/24),
((\xi_1,\xi_2,\xi_1),1/24),((\xi_2,\xi_1,\xi_1),1/24)$ \\
\hline
\end{tabular}
\end{center}
\end{table}

Table~\ref{sctable} gives the
values of relevant scalars for these test cases.
It should be noted that $\sigma_{\Q,p}=\tau_{\Q,p}=1$ for the $p=1$ cases
under the midpoint rule
because in these cases $S_{\Q,p}$ given by \eref{sqpdef}
turns out to be the identity matrix.

\begin{table}
\caption{Scalars relevant for test cases}
\label{sctable}
\begin{center}
\begin{tabular}{ccccccc}
\hline
Test & $d$ & $p$ & $q$ & $\sigma_{\Q,p}$ & $\tau_{\Q,p}$ & $M_\Q/m_q$ \\
\hline
$\TT_1$ & 2 & 1 & 1 & 1.00 & 1.00 & 1.0\\
$\TT_2$ & 2 & 2 & 3 & 5.26 & 0.83 & 1.0\\
$\TT_3,\TT_5$ & 3 & 1 & 1 & 1.00 & 1.00 & 1.0\\
$\TT_4,\TT_6$ & 3 & 2 & 4 & 6.47 & 0.63 & 1.0\\
\hline
\end{tabular}
\end{center}
\end{table}

Finally, we can tabulate element stiffness approximation bounds.
(As mentioned in Section~\ref{algo}, it suffices to measure element
stiffness approximation, which is much easier to compute than
global stiffness approximation.)  There are at least three
relevant quantities to tabulate: $\chi_1=\kappa(K_t,\bar{K}_t)$,
$\chi_2=\kappa(H)$, and $\chi_3=\hat\theta(\T)\kappa_1(\T)^2\kappa_2(\T)\cdot 
\frac{M_\Q\sigma_{\Q,p}^2}{m_\Q\tau_{\Q,p}^2}$, which is
the right-hand side of \eref{hbd}.  Here,
$K_t$ is the stiffness matrix of element $t$ ($t=1,\ldots,m$),
and $\bar{K}_t$ is the diagonally dominant preconditioner given
by $A_t^T\bar{D}_t A_t$.  It follows from
Theorem~\ref{thm:condno} and \eref{hbd} that $\chi_1\le\chi_2\le\chi_3$.

We generate six meshes and take the
max of $\chi_1,\chi_2,\chi_3$ over all elements of the mesh.  For
conductivity, we use $\theta\equiv 1$,
so that the factor $\hat\theta(\T)$ does not enter the bound.
The meshes are generated as follows.  In two dimensions,
we use the Triangle mesh generator \cite{Shewchuk} to generate
a mesh of an annulus, which is $\TT_1$ and consists of linear ($p=1$)
elements.
To generate $\TT_2$, a quadratic mesh,
we insert midpoint nodes in every element of $\TT_1$,
and for the edges on the boundaries, we moved the midpoint nodes
onto the boundary.

The remaining four meshes are three-dimensional.
For $\TT_3$, we used the QMG mesh generator \cite{MitchVava:QMG}
to generate a mesh of a unit ball in $\R^3$.  QMG's meshes are intended
for the $p=1$ case.  For $\TT_4$ we use midpoint
insertion on the $\TT_3$ mesh
to obtain a $p=2$ mesh.  Again, midpoint nodes of edges
whose endpoints were on the boundary
were displaced outward onto the boundary.  For $\TT_5$ and $\TT_6$, we
again mesh a unit ball (actually, one octant of the ball) using a
mesh generator tailored
for that domain only.  The mesh generator for $\TT_5$ first generates
a highly regular mesh of 
a unit tetrahedron,
and then it projects the mesh nodes 
radially outward toward one facet so that
an octant of the unit ball is covered.
Finally, $\TT_6$ is obtained from $\TT_5$ using midpoint insertion and
displacement at the boundary.

Then for every element in each mesh, we compute the three quantities
$\chi_1,\chi_2,\chi_3$. We have tabulated the maximum values in
Table~\ref{condtable}.  The undesirably large values of $\chi_1$ 
for $\TT_3$ and $\TT_4$
appears to be due primarily
to poorly shaped elements, i.e., large value of $\kappa_1(\T)$.
This is evident from comparing $\TT_3$ and $\TT_4$ against
$\TT_5$ and $\TT_6$, two meshes for the same domain but with
much better shaped tetrahedra.

The value of $\kappa_2(\T)$ was quite small in all tests; as
mentioned earlier, it is identically 1 for the $p=1$ cases. For
the $p=2$ cases, interior elements are still linear, and
elements adjacent to the boundary are fairly close to
linear.
It is interesting to note that
$\chi_1=\chi_2$ in all cases.  This implies that the invariant subspaces
of $H$ corresponding to its extremal eigenvalues meet the
range space of $\bar D^{1/2}A$.  We do not have a deeper explanation
for this observation.

\begin{table}
\caption{Approximation of the element stiffness matrix
by a diagonally dominant matrix.}
\label{condtable}
\begin{tabular}{rrrrrrrr}
\hline\\
Test & $m$ & $n$ & $\kappa_1(\T)$ & $\kappa_2(\T)$ & $\chi_1$ & $\chi_2$ & $\chi_3$ \\
\hline
$\TT_1$ &  234 & 143 & 4.1 & 1.0 & 16.6 & 16.6 & 16.6 \\
$\TT_2$ &  234 & 520 & 4.7 & 1.3 & 185.1 & 185.1 & 1117.2 \\
$\TT_3$ & 9078 & 1913 & 59.6 & 1.0 & 3549.5 & 3549.5 & 3549.5 \\
$\TT_4$ & 9078 & 13608 & 90.2 & 2.7 & $5.96\cdot 10^4$
& $5.96\cdot 10^4$ & $2.32\cdot 10^6$ \\
$\TT_5$ & 729 & 220 & 5.0 & 1.0  & 24.9 & 24.9 & 24.9 \\
$\TT_6$ & 729 & 1330 & 5.0 & 1.2 & 396.6 & 396.6 & 2644.5 \\
\hline
\end{tabular}
\end{table}

Although an iterative method was not tested, some conclusions can
still be drawn from Table~\ref{condtable}.  For example, if our
approximation were used on $\TT_6$ in an iterative setting, then the
slowdown would be at most a factor of $\sqrt{\kappa(\bar K_t,K_t)}$, i.e., a slowdown
of at most a factor of 20.  In practice, the method seems to work better
than that according to Avron et al.

\section{Open questions}
This work is the first to extend support-tree methods, which previously
have been shown to be good preconditioners for diagonally dominant
matrices with negative off-diagonal entries, to the class of finite element
matrices.  We have shown that the scope of the method includes
a standard scalar elliptic boundary value problem, but perhaps
the scope of finite element problems that can be tackled with
this method could be expanded further.

One generalization would be the class of problems $\nabla\cdot
(\Theta(x)\nabla u)=-f$, where $\Theta(x)$ is a spatially varying
$d\times d$ symmetric positive definite matrix.  This generalization
would present problems for our current analysis in the case that $\Theta(x)$
is highly ill-conditioned.  It would still be straightforward to write
$K=A^TJ^TDJA$ where $D$ is now
block diagonal, but our analysis of the introduction of $\bar J$ would
run into trouble because the $dq\times dq$
diagonal blocks of $D$ are no longer individually well conditioned.

It would also be interesting to tackle vector
problems such as linear elasticity or 
Stokes' flow, or higher-order
equations like the biharmonic equation.
It seems likely that our techniques can extend to
at least some of these problems since they all have a symmetric
positive definite weak form.
On the other hand,
a recent result \cite{ChenToledo} suggests that
the extension of our results to linear elasticity will not be
straightforward because of the higher nullity
of element stiffness matrices in the case of linear elasticity.
A further generalization would be to unsymmetric problems
like the convection-diffusion equation.  The latter class
of problems would require
substantial rethinking of the whole approach since condition number reduction,
which is very relevant for
the application of conjugate gradients to symmetric
positive definite systems, is less relevant to the application of GMRES
to unsymmetric systems.

Our analysis is based on condition numbers (support numbers). 
One drawback of this approach is that the convergence and work estimates 
may be too pessimistic. For instance,
the condition number of the preconditioned linear systems depends on
$\kappa_1(\T)$, the worst aspect ratio of any element in the mesh. 
If there is only one poorly shaped element in the mesh,
we expect iterative solvers will only take a few extra iterations since
changing a single element implies a low-rank correction to
the assembled stiffness matrix. Any analysis based on condition numbers will
be unable to capture this effect. A related open issue is whether we can exploit
recent work in mesh quality metrics \cite{Knupp} to show that ``good
meshes'' both have small error in the FEM approximation and also produce
linear systems that can be well approximated by diagonally dominant systems.

Another point to make about our method is that, although the condition
number of the preconditioned system has an upper bound independent of
$$R_\theta=\max_{x\in\Omega} \theta(x)/\min_{x\in\Omega}\theta(x),$$
there will
still be a loss of significant digits due to roundoff
error when using our method
in the case that $R_\theta$ is large.  This is because
the system matrix and the
preconditioner separately are ill-conditioned operators.  
A special case of \eref{bvp} in which
$R_\theta$
is extremely large and in which
$f=0$ was considered in \cite{Vava:fe}.
That paper proposed a method based on Gaussian elimination
in which the loss of significant
digits is avoided.  Some of the ideas behind \cite{Vava:fe}
were also extended to solution via conjugate gradient using
support preconditioners by \cite{HowleVava2}.
The methodology in  \cite{HowleVava2}, however, was for node-arc
incidence matrices and for a particular kind of support preconditioner
called a {\em support tree} \cite{Gre96},
and it is not clear whether that method would apply
to the present setting.

Finally, there is greater effort to couple automatic finite element
mesh generation to detailed analysis of the solver and its error bounds.
Our approximation result adds another ingredient to this mix; it would
be interesting if an automatic finite element mesh generator could be tailored
to reduce specifically the quantity $\kappa(K,\bar{K})$.

\section*{Acknowledgments}
The authors wish to thank the four anonymous referees of this
paper for their very helpful comments.

\bibliographystyle{plain}
\bibliography{../../Bibfiles/support,../../Bibfiles/support2}

\begin{thebibliography}{10}

\bibitem{Argyris}
J.~H. Argyris.
\newblock The natural factor formulation of the stiffness for the matrix
  displacement method.
\newblock {\em Computer Methods in Applied Mechanics and Engineering},
  5:97--119, 1975.

\bibitem{Avron}
H.~Avron, D.~Chen, G.~Shklarski, and S.~Toledo.
\newblock Combinatorial preconditioners for scalar elliptic finite-element
  problems.
\newblock Preprint, 2006.

\bibitem{Axelsson}
O.~Axelsson and V.~A. Barker.
\newblock {\em Finite element solution of boundary value problems: theory and
  computation}.
\newblock SIAM Press, Philadelphia, 2001.

\bibitem{BGHNT}
M.~Bern, J.~R. Gilbert, B.~Hendrickson, N.~Nguyen, and S.~Toledo.
\newblock Support-graph preconditioners.
\newblock {\em SIAM J. Matrix Anal. App.}, 27:930--951, 2006.

\bibitem{BH04}
E.~G. Boman and B.~Hendrickson.
\newblock Support theory for preconditioning.
\newblock {\em {SIAM} J. on Matrix Anal. and Appl.}, 25(3):694--717, 2004.
\newblock (Published electronically on 17 Dec 2003.).

\bibitem{CT03}
D.~Chen and S.~Toledo.
\newblock {V}aidya's preconditioners: Implementation and experimental study.
\newblock {\em {ETNA}}, 16, 2003.
\newblock Available from \url{http://etna.mcs.kent.edu}.

\bibitem{ChenToledo}
D.~Chen and S.~Toledo.
\newblock Combinatorial characterization of the null spaces of symmetric
  {H}-matrices.
\newblock {\em Linear Algebra and its Applications}, 392:71--90, 2004.

\bibitem{Cools}
R.~Cools.
\newblock An encyclopaedia of cubature formulas.
\newblock {\em J. Complexity}, 19:445--453, 2003.

\bibitem{EEST}
Michael Elkin, Yuval Emek, Daniel~A. Spielman, and Shang-Hua Teng.
\newblock Lower-stretch spanning trees.
\newblock In {\em Proceedings of the 37th annual ACM symposium on Theory of
  computing}, pages 494--503, New York, NY, USA, 2005. Association for
  Computing Machinery.

\bibitem{GVL}
G.~Golub and C.~F. Van~Loan.
\newblock {\em Matrix Computations, 3nd Edition}.
\newblock Johns Hopkins University Press, Baltimore, 1996.

\bibitem{Gre96}
K.~Gremban.
\newblock {\em Combinatorial Preconditioners for Sparse, Symmetric, Diagonally
  Dominant Linear Systems}.
\newblock PhD thesis, School of Computer Science, Carnegie-Mellon University,
  1996.
\newblock Available as Tech. Report CMU-CS-96-123.

\bibitem{gupta}
R.~Gupta.
\newblock Support graph preconditioners for sparse linear systems.
\newblock Master's thesis, Texas A\&M University, 2004.
\newblock Advised by V. Sarin. Available online at
  \url{https://txspace.tamu.edu/bitstream/1969.1/1353/1/etd-tamu-2004C-2-CPSC-%
Gupta.pdf}.

\bibitem{Gus94}
I.~Gustafsson.
\newblock An incomplete factorization preconditioning method based on
  modification of element matrices.
\newblock {\em BIT}, 36:86--100, 1996.

\bibitem{HowleVava2}
V.~E. Howle and S.~A. Vavasis.
\newblock An iterative method for solving complex-symmetric systems arising in
  electrical power modeling.
\newblock {\em SIAM J. Matrix Analysis App.}, 26:1150--1178, 2005.

\bibitem{Johnson}
C.~Johnson.
\newblock {\em Numerical Solution of Partial Differential Equations by the
  Finite Element Method}.
\newblock Cambridge University Press, 1987.

\bibitem{Knupp}
P.~Knupp.
\newblock Algebraic mesh quality metrics.
\newblock {\em SIAM J. Sci. Comput.}, 23:193--218, 2001.

\bibitem{MitchVava:QMG}
S.~A. Mitchell and S.~A. Vavasis.
\newblock Quality mesh generation in higher dimensions.
\newblock {\em SIAM J. Computing}, 29:1334--1370, 2000.

\bibitem{PhilMil}
T.~Phillips and G.~Miller.
\newblock Private communication.
\newblock 2005.

\bibitem{Shewchuk}
J.~R. Shewchuk.
\newblock Triangle: Engineering a {2D} quality mesh generator and {D}elaunay
  triangulator.
\newblock In Ming~C. Lin and Dinesh Manocha, editors, {\em Applied
  Computational Geometry: Towards Geometric Engineering}, volume 1148 of {\em
  Lecture Notes in Computer Science}, pages 203--222. Springer-Verlag, May
  1996.
\newblock From the First ACM Workshop on Applied Computational Geometry.

\bibitem{ST04}
Daniel~A. Spielman and Shang-Hua Teng.
\newblock Nearly-linear time algorithms for graph partitioning, graph
  sparsification, and solving linear systems.
\newblock In {\em Proceedings of the 36th annual ACM symposium on Theory of
  computing}, pages 81--90, New York, NY, USA, 2004. ACM Press.

\bibitem{Vava:bezier}
S.~Vavasis.
\newblock A {B}ernstein-{B}\'ezier sufficient condition for invertibility of
  polynomial mappings.
\newblock Archived by \url{http://arxiv.org/abs/cs.NA/0308021}, 2003.

\bibitem{Vava:fe}
S.~A. Vavasis.
\newblock Stable finite elements for problems with wild coefficients.
\newblock {\em SIAM J. Numer. Anal.}, 33:890--916, 1996.

\bibitem{WangSarin}
M.~Wang and V.~Sarin.
\newblock Parallel support graph preconditions.
\newblock In {\em High Performance Computing -- HiPC 2006}, volume 4297 of {\em
  Lecture Notes in Computer Science}, pages 387--398. Springer, Berlin, 2006.

\end{thebibliography}
\end{document}